\documentclass[12pt]{article}
\usepackage[T2B]{fontenc}
\usepackage{indentfirst}
\usepackage[pdftex, colorlinks, urlcolor=blue, pdfborder={0 0 0 [0 0]}]{hyperref}

\usepackage{amsfonts,amsmath,amssymb,theorem}
\usepackage{array}
\usepackage{esint}
\usepackage{graphicx}
\newtheorem{thm}{Theorem}[section]
\newtheorem{lem}{Lemma}[section]

\theoremstyle{remark}
\newtheorem{rem}{Remark}[section]
\theoremstyle{definition}

\title{Dynamics of random chains of finite
size with an infinite number of elements in $ {\mathbb R}^{2}$\footnote{%
\fbox{Published in \textcolor{blue}{\textit{Theory of Stochastic processes. 2010 vol.16 (32), no. 2.  Pp. 58-68}} }} }

\author{Elena V. Karachanskaya (Chalykh)\\
\small{\textit{Pacific National University,
Khabarovsk, Russia}}\footnote{%
Electronic address: Karachanckaya@mail.khstu.ru}}

\date{}

\begin{document}
\maketitle



\begin{abstract}
This article studies the dynamics of a finite chain with infinite
components. The equation which permits us to find the probability
distribution of the chain length is constructed and analysed. This
research is a continuation of paper \cite{1}.
\end{abstract}

\section{Formulation of the problem}

In Feller's book \cite{2} the problem of the length of a random
chain is considered, this chain is described in the following way:
the number of the elements is equal to $n$, the length of all its
elements is equal to one, the angle of one component with respect
to the previous is always the same up to a sign (the probability
of each angle is equal to 1/2), the distance between the end
points of the chain (length of the chain) is defined by means of
the average square length
$${\bf
M}[L_{n}^{2}]=n\,\displaystyle\frac{1+ \cos \alpha}{1-\cos
\alpha}-2\cos \alpha\,\displaystyle\frac{1- \cos^{n} \alpha}{(1-
\cos \alpha)^{2}}.$$

 We will consider the following chain: the
length of the chain is finite, the number of the components is
infinite, the length of each component is a random variable, the
angle of each component with respect to the previous one is also
random.
\begin{center}
\includegraphics[scale=0.5]{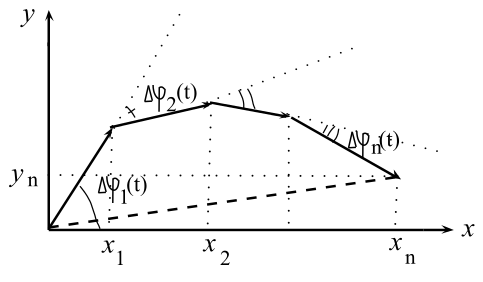}
\end{center}
The physical model can be a rope in a medium of Brownian
particles, as length of the chain we can understand the modulus of
the vector joining the starting point and the end point of the
chain.

Let $l \in [0,L] $ a parameter, $L$ a constant, $ l_{1},\,
l_{2},\, \ldots $ the values of the parameter, \,$
l_{1}<l_{2}<~\ldots~ \le~L $, \ $ \Delta = {L}/{N} $, \,
$l_{j}=j\cdot \Delta$. We will consider the model of the chain
described by the following system of equations:
\begin{equation}\label{1}
\begin{array}{c}\displaystyle
x_{N}(t) = \sum_{s=1} ^{N} a(l_{s})\Delta \cdot\cos{\varphi_{s}
(t)} , \ \ \ \ \ \displaystyle y_{N}(t) = \sum_{s=1} ^{N} a(l_{s}
) \Delta\cdot \sin{\varphi_{s} (t)},
\end{array}
\end{equation}
where $ a(l_{s}), \varphi_{s} (t) $ in general are random
processes, $ a(l)>0 $, $\varphi_{s} (t)$ is the angle between the
$s$-th component with respect to the previous one, $\varphi_{1}
(t)$ is the angle of the first component of the chain with respect
to the positive direction of the $x$-axis.

If we denote by $ |x_{n}(t) |^{2} + | y_{n}(t) |^{2} $ the length
of the chain consisting of the $n$ elements, then the length
$\triangle (l) $ of the real component of the chain is expressed
by the variable
$$ \triangle(l)
= a(l) \Delta, \ \ a(l)> 0, \ \
\displaystyle\int_{o}^{L}a(l)dl={\mathcal L}.
$$
Models of type \eqref{1} describe the distribution of the length $
{\mathcal L}(t) $ of the chain for the case where the following
inequality is satisfied:
$$ {\mathcal L}^{2}(t) = |x_{N}(t) |^{2} + | y_{N}(t) |^{2} \le const.
$$
From the point of view of the representation of the phenomenon of
the turbulent diffusion, the model \eqref{1} can be useful for
some generalizations of the passive displacement under the action
of vortices of different size \cite{3}.

Let be $n<N$ (that is we will consider not the whole chain but a
part of it), $N\to \infty$. Since the coordinates of the initial
point and the end point of each component depend on time $t$ and
from the parameter $l$, then in the model \eqref{1} we introduce
some changes.
\begin{equation}\label{2}
\begin{array}{c} \displaystyle
x_{n}(l;t) = \sum_{s=1} ^{n} a(l_{s})\Delta \cdot\cos{\varphi_{s}
(t)} , \ \ \ \ \ \displaystyle y_{n}(l;t) = \sum_{s=1} ^{n}
a(l_{s} ) \Delta\cdot \sin{\varphi_{s} (t)}.
\end{array}
\end{equation}

In this way, the random field $ \{ x_{n}(l;t); y_{n}(l;t) \} $ is
a dynamical stochastic process. We will study its limit behavior
for $n \to \infty$.

\section{Assumptions on the model}

In order to obtain the coefficient of the limit equation in
analytical form we shall restrict ourself to the model satisfying
the following assumptions:
\begin{equation}\label{3}
\begin{array}{c}
  a(l)>0 , \ \ \ l \in [0,L],\\
  \varphi_{s}(t) = \displaystyle\sum\limits_{k=1}^{s} \eta (l_{k}; t) \Delta(w(l_{k})),
\ \ \ t \in [0,T], \\
\eta (l_{k}; t) = \displaystyle\int_{0}^{t} \sigma (l_{k}; \tau )
\, d
w_{k} (\tau ),\\
\end{array}%
\end{equation}
where $\triangle(w(l_{k})), $ $ \triangle (w_{k}(\tau )) $ are
independent among themselves and for different $s$ and $\tau$ are
anticipating increments of the corresponding Wiener processes
defined on the product of independent probability spaces
$$
\{ \Omega_{1},\, \Im_{l}, \, P_{1} \} \, \times \, \{ \Omega_{2},
\, \Im_{t}(n), \, P_{2} \},
$$
where  $ \Im_{l} , \ $ and $ \Im_{t}(n) \ $ are the corresponding
flows of sigma algebras generated by the processes $ w(l) $ and $
w(t) \in {\mathbb R}^{n}  $; the functions $ a(l)\in {\mathbb
C}^{1}_{[0,L]}$  and $ \sigma (l; t) \in {\mathbb
C}^{2}_{[0,L]\times [0,T]} $ are deterministic functions depending
on $l$ and $t$, $\eta (l_{s}; t)$ is the intensity of the angle.

We have therefore,
\begin{equation}\label{4}
\begin{array}{c}
   \displaystyle
x_{n}(l;t) = \sum_{s=1} ^{n} a(l_{s})\Delta \cdot
\cos\left[{\sum\limits_{k=1}^{s} \left(\int_{0}^{t} \sigma (l_{k};
\tau ) \, d
w_{k} (\tau )\right)\triangle (w(l_{k}))}\right]\\
  \displaystyle
y_{n}(l;t) = \sum_{s=1} ^{n} a(l_{s})\Delta \cdot
\sin\left[{\sum\limits_{k=1}^{s} \left(\int_{0}^{t} \sigma (l_{k};
\tau ) \, d w_{k} (\tau )\right)\triangle (w(l_{k}))}\right]
  \end{array}
\end{equation}

Under the condition of bounded length of the chain for the random
function $\varphi_{s}(t)$ can be defined the limit for $n \to
\infty$. In this context the variable $l$ appears as a parameter.

\section{Transition to auxiliary processes}

Let us transform (\ref{2}) by means of the Euler representation:
\begin{equation*}
\begin{array}{c}
 x_{n}(l;t)=\displaystyle\sum_{s=1} ^{n} a(l_{s})\Delta
\cdot\cos{\varphi_{s} (t)}=\sum_{s=1} ^{n} a(l_{s})\Delta \cdot
\displaystyle\frac{\exp\{i\varphi_{s} (t)\}+\exp\{-i\varphi_{s}
(t)\}}{2}=
\\
=\displaystyle\frac{1}{2}\sum_{s=1} ^{n} a(l_{s})\Delta \cdot
\exp\{i\varphi_{s} (t)\} +\displaystyle\frac{1}{2}\sum_{s=1} ^{n}
a(l_{s})\Delta \cdot \exp\{-i\varphi_{s} (t)\}, \\
y_{n}(l;t)=\displaystyle\sum_{s=1} ^{n} a(l_{s})\Delta
\cdot\sin{\varphi_{s} (t)}=\sum_{s=1} ^{n} a(l_{s})\Delta \cdot
\displaystyle\frac{\exp\{i\varphi_{s} (t)\}-\exp\{-i\varphi_{s}
(t)\}}{2i}=
\\
=\displaystyle\frac{1}{2i}\sum_{s=1} ^{n} a(l_{s})\Delta \cdot
\exp\{i\varphi_{s} (t)\} -\displaystyle\frac{1}{2i}\sum_{s=1} ^{n}
a(l_{s})\Delta \cdot \exp\{-i\varphi_{s} (t)\},
\end{array}
\end{equation*}
we now introduce the auxiliary process
\begin{equation*}
\begin{array}{c}
  z_{1}(s;t) = \exp \{ - i \displaystyle\sum_{j=1}^{s} \, \triangle w(l_{j})
\displaystyle\int_{0}^{t} \sigma (l_{j}; \tau ) \, d w_{j} ( \tau
) \} ,
\\
z_{n,1}(l,t) = \displaystyle\sum_{s=1}^{n} a(l_{s})\Delta \cdot
\exp\{-i\varphi_{s} (t)\} =\sum_{s=1}^{n}  a(l_{s})  \Delta \cdot
z_{1}(s;t), \ \ \Delta = O(n^{-1}).
\end{array}
\end{equation*}
By using the Euler representation we rewrite the
process $ \{ x_{n}(l;t); $ $ y_{n}(l;t) \} $  in the following
form
\begin{equation*}%\label{5}
x_{n}(l;t) = \frac{1}{2} (z_{n,1}(l,t) + z^{*}_{n,1}(l,t)),  \ \ \
y_{n}(l;t) = \frac{i}{2} (z_{n,1}(l,t) - z^{*}_{n,1}(l,t)).
\end{equation*}
For the construction of the characteristic function of the random
field $ \{ x_{n}(l;t); $ $ y_{n}(l;t) \} $ we define the form of
the function $\exp\{i(\alpha x_{n}(l;t) +\beta y_{n}(l;t))\}$:
\begin{equation*}
\begin{array}{c}
\displaystyle\exp\{i(\alpha x_{n}(l;t) +\beta y_{n}(l;t))\}=\exp{
\left\{ i \alpha \frac{z_{n,1}(l;t) + z_{n,1}^{*}(l;t)}{2} - \beta
\frac{z_{n,1}(l;t) - z_{n,1}^{*}(l;t)}{2} \right\} } =
\\
= \displaystyle\sum \limits _{m,r=1}^{\infty } \frac{(i\alpha -
\beta)^{m} (i\alpha + \beta)^{r}}{2^{m+r} m! r!} z_{n,1}^{m}(l;t)
z_{n,1}^{*^{r}}(l;t).
\end{array}
\end{equation*}
In consequence the analysis of the process  $ \{ x_{n}(l; t);
y_{n}(l;t) \} $  leads to the study of the process $
z_{n,1}^{m}(l;t) z_{n,1}^{*^{r}}(l;t)$.

Since the summation and the integration operations have the same
properties we replaced (in symbolic form, when $\Delta \to 0$ and
this corresponds to $n \to \infty$) the process
\begin{equation*}
z_{n,1}(l,t) = \displaystyle\sum_{s=1}^{n}  a(l_{s})  \Delta \cdot
\exp \left\{ - i \sum_{j=1}^{s}\left(  \int\limits _{0}^{t} \sigma
(l_{j}; \tau ) \, d w_{j} ( \tau )\right)\, \triangle w(l_{j})
\right\}
\end{equation*}
by the process
\begin{equation*}%\label{6}
\begin{array}{c}
\displaystyle z_{,1} = \sum \limits_{s=1}^{n} a(l_{s}) \Delta
\cdot\exp \left\{ -\sum \limits_{j=1}^{s} \eta (l_{s},t) \triangle
w(l_{s})\right\}=
%\\ =
\displaystyle\fint\limits_{0}^{l}a(u)\exp \left\{-i
\fint\limits_{0}^{u}\eta (\theta,t) d w(\theta) \right\}du,
\end{array}
\end{equation*}
where $\eta(u,t)=\displaystyle\int_{0}^{t}\sigma (u,\tau)d\tau$.
We do not loose any generality in the analysis with this
assumption and in the sequel we shall use the symbol
$\displaystyle\fint$ instead of
$\sum$.\footnote{%
This symbol does not concern known designations. It is a label
only.}

\section{Degree transformation}

By considering the continuity of the process $ z_{n,1}(l;t) $ and,
consequently, of the process $ z_{\,,1}(l;t) $ with respect to
both variables $ l$ and $t$, we produce the degree transformation:
\begin{equation}\label{7}
\begin{array}{c} \displaystyle z_{,1}^{m}(l;t) =
\left[ \fint_{0}^{l} a(u) \exp \left \{ -i \fint_{0}^{u} \eta
(\theta;t) \, dw(\theta ) \right\} \, du \right]^{m}=
\\
\displaystyle =  m! \fint_{0}^{l} a(u_{1}) \, du_{1} \, \exp \left
\{ -i\, m \fint_{0}^{u_{1}} \eta (\theta_{1}; t) \, dw(\theta_{1}
) \right\}  \times
\\
\displaystyle \times \fint_{u_{1}}^{l} a(u_{2}) \, du_{1} \,\exp
\left\{ -i(m-1) \fint_{u_{1}}^{u_{2}} \eta (\theta_{2};t) \,
dw(\theta_{2} ) \right\} \times \, \ldots \, \times
\\
%\displaystyle \times \fint_{u_{m-2}}^{l} a(u_{m-1}) \, du_{m-1}
%\,\exp \left \{ -2i \fint_{u_{m-2}}^{u_{m-1}} \eta
%(\theta_{m-1};t) \, dw(\theta_{m-1} ) \right\} \times
%\\
\displaystyle \times \fint_{u_{m-1}}^{l} a(u_{m}) \, du_{m} \,\exp
\left \{ -i \fint_{u_{m-1}}^{u_{m}} \eta (\theta_{m};t) \,
dw(\theta_{m} ) \right\}
\end{array}
\end{equation}
where $ 0 < u_{1} < \, \ldots \, < u_{m} <l$. In different
intervals we have different $dw(\theta)$ for each time instant
$t$.

\section{Determination of moments}

Since the process $ z^{m}_{,1}(l;t) z^{*{r}}_{,1}(l;t) $ depends
from two variables, for the calculation of the mean $ {\bf
M} \, [ z^{m}_{,1}(l;t) z^{*{r}}_{,1}(l;t) ] $ it is
necessary to carry out the averaging for $t$ constant (on the
space $\Omega_1$), and then $\Omega_l$.

\subsection{Averaging with respect to $l$}

On the disjoint intervals $ [u_{i}, u_{i+1} )$ for all $ i $ , for
all $i$ the processes
$$ f(u_{i}, u_{i+1}) =
\fint_{u_{i}}^{u_{i+1}} \eta (\theta_{i+1};t) \, dw(\theta_{i+1} )
 $$
are independent by construction (since $ dw(\theta )$ are Wiener
processes,  $\eta (\theta;t) $ for fixed $t$ is a non-random
function depending on $\theta$). Because of this, since $a(l)$ is
also a non-random function, then the mathematical mean of each
factor is defined in the following way:
\begin{equation}\label{8}
\begin{array}{c}
\displaystyle{\bf M}_{t} \, [ z_{,1}^{m}(l;t) ] = m!
\fint\limits_{0}^{l} a(u_{1}) \, du_{1} \,{\bf M}_{t} \left[ \exp
\left \{ -i\, m \fint_{0}^{u_{1}} \eta (\theta_{1}; t) \,
dw(\theta_{1} ) \right\} \right] \times
\\
\displaystyle \times \fint_{u_{1}}^{l} a(u_{2}) \, du_{2} \,{\bf
M}_{t} \left[ \exp \left\{ -i(m-1) \fint_{u_{1}}^{u_{2}} \eta
(\theta_{2};t) \, dw(\theta_{2} ) \right\} \right] \times \,
\ldots \, \times
\\
%\displaystyle \times \fint_{u_{m-2}}^{l} a(u_{m-1}) \, du_{m-1}
%\,{\bf M}_{t} \left[ \exp \left\{ -2i \fint_{u_{m-2}}^{u_{m-1}}
%\eta (\theta_{m-1};t) \, dw(\theta_{m-1} ) \right\} \right]
%\times
%\\
\displaystyle \times \fint_{u_{m-1}}^{l} a(u_{m}) \, du_{m} \,{\bf
M}_{t} \left[ \exp \left\{ -i \fint_{u_{m-1}}^{u_{m}} \eta
(\theta_{m};t) \, dw(\theta_{m} ) \right\} \right] .
\end{array}
\end{equation}
In this way it is necessary to find the mean of the expression of
the following type:
$$
\exp{ \left\{ -i (m-j) \fint_{u_{j}}^{u_{j+1}} \eta (\theta;t) \,
dw(\theta) \right\} } .
$$
\begin{lem}\label{lm1}
{\it The following equality holds
\begin{equation}\label{9}
\displaystyle {\bf M}_{t} \, \left[ \exp{ \left\{ \alpha
\fint_{a}^{b} \eta (u;t)\, dw(u) \right\} } \right] = \exp{
\left\{ \frac{1}{2} \alpha^{2} \fint_{a}^{b} \eta^{2} (u;t)\, du
\right\} } .
\end{equation}}
\end{lem}
{\it Proof.} Let us denote:
\begin{equation}\label{10}
q(a,b;t) = \fint_{a}^{b} \eta (u;t)\, dw(u)
\end{equation}
and differentiate $ q(a,b;t) $ with respect to the upper limit
$b$. As a result we obtain:
$$
d_{b} \, q(a,b;t) = \eta (b;t) \, dw(b).
$$
Therefore, by Ito formula, the stochastic differential with
respect to $b$ of the expression
$$ \exp{ \left\{ \alpha \fint_{a}^{b} \eta (u;t)\, dw(u) \right\} }=\exp{\{ \alpha
q(a,b;t) \} } $$
is equal to
$$
d_{b} \, \exp{\{ \alpha q(a,b;t) \} } = \exp{\{ \alpha q(b;t)\} }
\alpha \eta (b;t)  \, dw(b)+ \frac{1}{2} \alpha^{2} \eta^{2} (b;t)
\exp{ \{ \alpha q(a,b;t)\} } db.
$$
We compute the average with respect to $l$ of the obtained
expression:
\begin{equation}\label{11}
d_{b} \, {\bf M}_{t} \, \left[ \exp{\{ \alpha q(a,b;t) \} }
\right]= \frac{1}{2} \alpha^{2} \eta^{2} (b;t) {\bf M}_{t} \,
\left[ \exp{ \{ \alpha q(a,b;t)\} } \right]db.
\end{equation}
Let us denote:
\begin{equation}\label{12}
I_{1}(a,b;t) =  {\bf M}_{t} \, \left[ \exp{\{ \alpha q(a,b;t)\} }.
\right]
\end{equation}
Let $ \eta (b;t) $ be independent from the stochastic process $
w(u) $. In view of \eqref{11}  we obtain the differential equation
$$
\displaystyle\frac{d I_{1}(a,b;t)}{db}=\frac{1}{2} \alpha^{2}
\eta^{2} (b;t) I_{1}(a,b;t)
$$
and its solution
$$
I_{1}(a,b;t) = \exp{ \left\{ \frac{1}{2} \alpha^{2} \fint_{a}^{b}
\eta^{2} (u;t)\, du \right\} },
$$
satisfies the initial condition $I_{1}(a,a;t)=1$. In view of
\eqref{10}, \eqref{12}  the statement of the lemma is proved. \ \
\ $\square$

As a consequence of Lemma \ref{lm1}, the mathematical mean
\eqref{8} takes the form (for $t $ constant):
$$
\begin{array}{c}
 {\bf M}_{t} \, [ z_{,1}^{m}(l;t) ] = m! \displaystyle\fint_{0}^{l}
a(u_{1}) \, du_{1} \, \exp{ \left\{ - \frac{m^{2}}{2}
\fint_{0}^{u_{1}} \eta^{2} (\theta; t) \, d\theta \right\} }
\times
\\
 \times \displaystyle \fint_{u_{1}}^{l} a(u_{2}) \, du_{2} \, \exp{
\left\{ - \frac{(m-1)^{2}}{2} \fint_{u_{1}}^{u_{2}} \eta^{2}
(\theta;t) \, d\theta \right\} } \times \, \ldots \, \times
\\
\times
 \displaystyle\fint_{u_{m-1}}^{l} a(u_{m}) \, du_{m} \, \exp{
\left\{ - \frac{1}{2} \fint_{u_{m-1}}^{u_{m}} \eta^{2} (\theta;t)
\, d\theta \right\} }  .
\end{array}
$$

\subsection{Averaging with respect to $t$}

Now we make the averaging of the process $ z_{,1}^{m}(l;t) $ on
the space $\Omega_2$. Since $ w_{s}(t)$ for all $s$ are
independent Wiener processes, then $ \eta^{2} (l_{s}; t) $ are
independent, the average of the product is therefore equal to the
product of the means. As a result we obtain:
\begin{equation*}
\begin{array}{c}
 \displaystyle{\bf M} \, \left[ {\bf M}_{t} \, [ z_{,1}^{m}(l;t) ] \right] =
%$$ %\vskip0.5mm $$=
m! \fint_{0}^{l} a(u_{1}) \, du_{1} \, {\bf M}\, \left[ \exp{
\left\{ - \frac{m^{2}}{2} \fint_{0}^{u_{1}} \eta^{2} (\theta; t)
\, d\theta \right\} } \right] \times
%\\
%\times \fint_{u_{1}}^{l} a(u_{2}) \, du_{2} \, {\bf M} \, \left[
%\exp{ \left\{ - \frac{(m-1)^{2}}{2} \fint_{u_{1}}^{u_{2}} \eta^{2}
%(\theta;t) \, d\theta \right\} } \right] \times \, \ldots \,
%\times
\\
\displaystyle\times \, \ldots \,\fint_{u_{m-1}}^{l} a(u_{m}) \,
du_{m} \, {\bf M} \, \left[ \exp{ \left\{ - \frac{1}{2}
\fint_{u_{m-1}}^{u_{m}} \eta^{2} (\theta;t) \, d\theta \right\} }
\right].
\end{array}
\end{equation*}
\begin{lem}\label{lm2}
{\it The following relationship holds
$$
{\bf M} \, \left[ \exp{ \left\{ - \frac{\alpha^{2}}{2}
\fint_{a}^{b} \eta^{2} (u;t)\, du \right\} } \right] = \exp{
\left\{ - \frac{\alpha^{2}}{4} \fint_{a}^{b} \left( \int_{0}^{t}
\sigma^{2}(u; \tau ) \, d \tau \right) \, du \right\} } ,
$$
where
$$
\eta (u;t) = \int_{0}^{t} \sigma (u; \tau ) \, d w(\tau ),
$$
$ \sigma (u;t) $ is a non-random function.}
\end{lem}

{\it Proof.} We shall use the following representation:
$$
\fint_{a}^{b} \eta^{2}(u; t) \, du = \frac{b-a}{N} \sum
\limits_{k=1}^{N} \eta^{2}(u_{k};t),
$$
this is possible in force of the model assumptions. The processes
$ \eta^{2}(u_{k};t)$ \ ($k=\overline{1,N} $) by definition are
independent for different values of $k$. We introduce the
notation:
$$ P_{k}(t)=\eta (u_{k};t) = \int_{0}^{t} \sigma
(u_{k}; \tau ) \, d w_{k}(\tau ),
$$
where $ \sigma (u_{k};t) $ is a non-random function depending on
$u_k$, $t$. We consider now two cases.

\textbf{A.} Let $ \sigma (u_{k};t) $ be constant. For the seek of
simplicity in the sequel we assume that $\sigma(u_{k};t) =1$ and
study the problem for the processes
\begin{equation}\label{13}
\eta (u_{k};t) = \int_{0}^{t}dw_{k}(\tau) = {\tilde P}_{k} (t) .
\end{equation}
By considering the representation of the integral in form of sums,
we carry out the transformation:
$$
\exp{ \left\{ - \frac{\alpha^{2}(b-a)}{2N} \eta^{2}(u_{k},t)
\right\}}=\exp{ \left\{ - \frac{\alpha^{2}(b-a)}{2N}
\left[\int_{0}^{t} dw_{k}(\tau)\right]^{2} \right\}}=\exp{ \left\{
- \frac{\alpha^{2}(b-a)}{2N} {\tilde P}_{k}^{2} (t) \right\}}.
$$
Therefore
$$ \exp{ \left\{ - \frac{\alpha^{2}(b-a)}{2N}\sum_{k=1}^{N}
{\tilde P}_{k}^{2} (t) \right\}}=\prod\limits_{k=1}^{N}\exp{
\left\{ - \frac{\alpha^{2}(b-a)}{2N} {\tilde P}_{k}^{2} (t)
\right\}}
 $$
We denote by:
$$
\exp{ \left\{ - \frac{\alpha^{2}(b-a)}{2N}\sum\limits_{k=1}^{N}
{\tilde P}_{k}^{2} (t) \right\}}=I_{N}(k,\alpha^{2}).
$$
Since
$$
{\bf M} \, \left[ \exp{ \left\{\alpha q(a,b;t) \right\} }
\right]={\bf M} \, \left[ \exp{ \left\{-\alpha q(a,b;t) \right\} }
\right],
$$
and the following relationships hold:
\begin{equation*}
\begin{array}{c}
{\bf M} \, \left[\prod\limits_{k=1}^{N} \exp{ \left\{ -
\displaystyle\frac{\alpha^{2}(b-a)}{2N} {\tilde P}_{k}^{2} (t)
\right\}}
%\right\}}
 \right]=\\
 =\displaystyle\prod\limits_{k=1}^{N}{\bf M} \, \left[\exp{ \left\{ -
\frac{\alpha^{2}(b-a)}{2N} {\tilde P}_{k}^{2} (t)
%\right\}}
\right\} } \right] \ \ {\mathop\rightarrow^{qm}}\ \ {\bf M} \,
\left[ \exp{ \left\{\alpha q(a,b;t) \right\} } \right].
\end{array}
\end{equation*}
Then we carry out the Ito differentiation:
$$
d {\tilde P}_{k} (t) = dw_{k} (t),
$$
\begin{equation}\label{14}
d {\tilde P}_{k}^{2} (t) = dt + 2w_{k} (t)\, dw_{k} (t)
\end{equation}
and in view of \eqref{13}  we have that:
\begin{equation*}
\begin{array}{c} d_{t}  \displaystyle\exp{ \left\{ - \frac{\alpha^{2}(b-a)}{2N}\,
{\tilde P}_{k}^{2} (t) \right\}} = - \exp{ \left\{ -
\displaystyle\frac{\alpha^{2}(b-a)}{2N}\, {\tilde P}_{k}^{2} (t)
\right\}} \frac{\alpha^{2}(b-a)}{2N}\,d {\tilde P}_{k}^{2} (t) +
\\
+ \displaystyle\frac{\alpha^{4}(b-a)^{2}}{2N^{2}} \,{\tilde
P}_{k}^{2} (t) \exp{ \left\{ -\frac{\alpha^{2}(b-a)}{2N} \,
{\tilde P}_{k}^{2} (t) \right\} } \, dt.
\end{array}
\end{equation*}
We introduce in the last differential the expression \eqref{14}:
\begin{equation*}
\begin{array}{c}
d_{t}  \displaystyle\exp{ \left\{ -
\displaystyle\frac{\alpha^{2}(b-a)}{2N}\, {\tilde P}_{k}^{2} (t)
\right\} } = - \exp{ \left\{ - \frac{\alpha^{2}(b-a)}{2N}\,
{\tilde P}_{k}^{2} (t) \right\} } \times
\\
\times \displaystyle\left[ \displaystyle\frac{\alpha^{2}(b-a)}{2N}
\, 2{\tilde P}_{k} (t) \, dw_{k} (t) - \frac{\alpha(b-a)}{2N} \,
\left( -1 + \frac{\alpha^{2}(b-a)}{N} {\tilde P}_{k}^{2} (t)
\right) \, dt \right] .
\end{array}
\end{equation*}
We calculate the mean for the last expression by denoting
\begin{equation*}
I_{2}(t;\alpha^{2}) = {\bf M} \left[ \exp{ \left\{ -
\displaystyle\frac{\alpha^{2}(b-a)}{2N}\, {\tilde P}_{k}^{2} (t)
\right\} } \right] .
\end{equation*}
We obtain the equation:
\begin{equation*}
dI_{2}(t; \alpha^{2}) = I_{2}(t; \alpha^{2})
\, \frac{ \alpha^{2}(b-a)}{2N} \, dt +
%$$ $$ +
\frac{ \alpha^{4}(b-a)^{2}}{2N^{2}} \, {\bf M} \left[ {\tilde
P}^{2}_{k} (t) \exp{ \left\{ - \frac{\alpha^{2}(b-a)}{2N}\,
{\tilde P}_{k}^{2} (t) \right\}} \right] \, dt .
\end{equation*}
By considering the differentiation with respect to $ \alpha^{2} $
of the expression
\begin{equation*}
\exp{ \left\{ - \frac{\alpha^{2}(b-a)}{2N}\, {\tilde P}_{k}^{2}
(t) \right\} } ,
\end{equation*}
the last equation can be represented as a partial differential
equation with constant coefficients:
\begin{equation}\label{15}
\frac{ dI_{2}(t; \alpha^{2})}{dt} = - \frac{
\alpha^{2}(b-a)}{2N}\, I_{2}(t; \alpha^{2}) -  \frac{
\alpha^{4}(b-a)}{N} \, \frac { \partial} {\partial \alpha^{2}} \,
I_{2}(t; \alpha^{2}) .
\end{equation}
The solution of this equation will be obtained by exploiting the
properties of the stochastic processes. With this purpose we
evaluate the mean of the function \linebreak $ \exp{ \left\{ -
\displaystyle\frac{\alpha^{2}(b-a)}{2N}\, {\tilde P}_{k}^{2} (t)
\right\} }$. By considering that the process is a Wiener process
(that is a Gaussian process) we have that:
$$
{\bf M} \left[ \exp{ \left\{ - \frac{\alpha^{2}(b-a)}{2N}\,
{\tilde P}_{k}^{2} (t) \right\}} \right] =
%$$ $$ =
\int_{-\infty}^{\infty} \exp{ \left\{ -
\frac{\alpha^{2}(b-a)}{2N}\, x^{2} \right\} } \times \frac{1}{
\sqrt{2 \pi \, t}} \exp{ \left\{ - \frac{x^{2}}{2 t} \right\} } \,
dx =
$$
$$
= \frac{1}{  \sqrt{2 \pi \, t}} \, \int_{-\infty}^{\infty} \exp{
\left\{ - \left( \frac{\alpha^{2}(b-a)}{2N}\, x^{2} +
\frac{x^{2}}{2 t} \right) \right\} } \, dx .
$$
Furthermore
$$
{\bf M} \left[ \exp{ \left\{ - \frac{\alpha^{2}(b-a)}{2N}\,
{\tilde P}_{k}^{2} (t) \right\}} \right] = \frac{1}{\sqrt{2 \pi \,
t}} \cdot \sqrt{2 \pi} \cdot \sqrt{
t\Bigl/\left(\frac{\alpha^{2}(b-a)}{2N} t + 1\right)\Bigl.} =
%$$ $$=
\left( \frac{\alpha^{2}(b-a)}{2N} t + 1 \right)^{- {1/2}}.
$$
In this way, the obtained expression $ \displaystyle I_{2}(t;
\alpha^{2}) = \left( \frac{\alpha^{2}(b-a)}{2N} t + 1 \right)^{-
{1/2}} $ is the solution of the differential equation \eqref{15}.
Besides this, in view of Lemma \ref{lm1}, we have that:
\begin{equation*}
\begin{array}{c}
  {\bf M}_{t} \, \displaystyle\left[ \exp{\{ -\alpha q(a,b;t)\}})
\right] = \exp{ \left\{ - \frac{\alpha^{2}}{2}\int_{a}^{b}
\eta^{2} (u,t) \right\}} = \exp{ \left\{ -
\frac{\alpha^{2}(b-a)}{2N}\sum\limits_{k=1}^{N} \eta^{2} (u_{k},t)
\right\}}=
\\
=\displaystyle\prod\limits_{k=1}^{N} \exp{ \left\{ -
\frac{\alpha^{2}(b-a)}{2N}\eta^{2} (u_{k},t) \right\}}.
\end{array}
\end{equation*}
Therefore
\begin{equation*}
\begin{array}{c}
 \displaystyle{\bf M}\left[ {\bf M}_{t} \, \left[ \exp{\{ -\alpha
q(a,b;t)\}}) \right] \right]\ \ {\mathop\rightarrow^{qm}}\ \
\prod\limits_{k=1}^{N} {\bf M}\left[\exp{ \left\{ -
\frac{\alpha^{2}(b-a)}{2N}{\tilde P}_{k}^{2} (t)
\right\}}\right]= \\
 =\displaystyle\prod\limits_{k=1}^{N}\left(
\frac{\alpha^{2}(b-a)}{2N} \, t + 1 \right)^{- 1/2}=\left(
\frac{\alpha^{2}(b-a)}{2N} \, t + 1 \right)^{- {N/2}}
\end{array}
\end{equation*}
Passing to the limit we obtain the complete averaging with respect
to both components:
\begin{equation*}
\displaystyle{\bf M} \, \left[ I_{1}(a,b;t) \right] =
%$$ $$ =
\lim \limits_{N \rightarrow\infty} \, \left(
\frac{\alpha^{2}(b-a)}{2N} \, t + 1 \right)^{- {N/2}} = \exp {
\left\{ - \frac{\alpha^{2}(b-a)t}{4} \right\} } .
\end{equation*}

\textbf{B.} Now we consider the case $ \sigma (u_{k};t) \not= 1$.
In this case for $ P_{k} (t) $, we obtain the expression:
$$
d_{t} \, P_{k} (t) = \sigma (u_{k};t) \, dw_{k} (t) .
$$
In this way we have that
$$
d_{t} P_{k}^{2} (t) = \frac {\sigma^{2} (u_{k};t)}{2} 2dt + 2
P_{k} (t) \, \sigma (u_{k};t) \, dw_{k} (t)
$$
and
\begin{equation*}
\begin{array}{c}
d_{t} \displaystyle \exp{ \left\{ - \frac{\alpha^{2}(b-a)}{2N}\,
P_{k}^{2} (t) \right\}} = - \exp{ \left\{ -
\frac{\alpha^{2}(b-a)}{2N}\, P_{k}^{2} (t) \right\} } \times
\\
\times \displaystyle\frac{\alpha^{2}(b-a)}{2N} \, \left[
\sigma^{2} (u_{k};t) \, dt + 2 P_{k} (t)\, \sigma (u_{k};t) \,
dw_{k} (t) \right] +
\\
+ \displaystyle\frac{\alpha^{4}(b-a)^{2}}{2N^{2}} \, P_{k}^{2} (t)
\, \sigma^{2} (u_{k};t) \cdot \exp{ \left\{
-\frac{\alpha^{2}(b-a)}{2N} \, P_{k}^{2} (t) \right\} } \, dt .
\end{array}
\end{equation*}
Therefore
\begin{equation*}
\begin{array}{c}
\displaystyle d_{t} I_{t} (k; \alpha^{2}) = \partial_{t} \, {\bf
M} \, \left[ \exp{ \left\{
-\displaystyle\frac{\alpha^{2}(b-a)}{2N} \, P_{k}^{2} (t) \right\}
} \right] =
\\
= -  \displaystyle{\bf M} \,\left[ \exp{ \left\{
-\displaystyle\frac{\alpha^{2}(b-a)}{2N} \, P_{k}^{2} (t) \right\}
} \frac{\alpha^{2}(b-a)}{2N} \, \sigma^{2} (u_{k};t) \right] \, dt
.
\end{array}
\end{equation*}
By exploiting the possibility of differentiation with respect to
the parameter $ \alpha^{2} $, we arrive at the following equation:
\begin{equation}\label{16}
\displaystyle \frac { \partial_{t} I_{t} (k; \alpha^{2}) }{
\partial t} = - \frac{\alpha^{2}(b-a)}{2N} \, \sigma^{2} (u_{k};t)
\, I_{t} (k;\alpha^{2}) -
%$$ \vskip0.5mm $$ -
\frac{\alpha^{4}(b-a)}{N} \, \sigma^{2} (u_{k};t) \, \frac {
\partial I_{t} (k; \alpha^{2}) }{ \partial \alpha^{2}} .
\end{equation}
We divide both members by $ \sigma^{2} (u_{k};t)$, and denote
$$
\theta (t) = \int_{0}^{t} \, \sigma^{2} (u_{k};\tau) \, d \tau
$$
and we pass to the auxiliary equation
\begin{equation}\label{17}
\displaystyle \frac { \partial I_{\theta} (k; \alpha^{2})
}{\partial \theta} = - \frac{\alpha^{2}(b-a)}{2N} \, I_{\theta}
(k;\alpha^{2}) - \frac{\alpha^{4}(b-a)}{N} \, \frac { \partial
I_{\theta} (k; \alpha^{2}) }{\partial \alpha^{2}} .
\end{equation}
Equation \eqref{17}  is a differential equation with constant
coefficients stochastically equivalent to equation \eqref{15} .
Therefore its solution has the form:
$$
I_{\theta} (k; \alpha^{2}) = \left( \frac{\alpha^{2}(b-a)}{2N}\,
\theta + 1 \right)^{- {1/2}}.
$$
and we obtain, the solution of the equation \eqref{16}:
$$
I_{t} (k; \alpha^{2}) = \left( \frac{\alpha^{2}(b-a)}{2N}\,
\int_{0}^{t} \, \sigma^{2} (u_{k};\tau) \, d \tau  + 1 \right)^{-
{1/2}}.
$$
As a consequence we have that
\begin{equation}\label{18}
\displaystyle{\bf M} \, \left[ I_{N}(k;\alpha^{2}) \right]
=\prod\limits_{k=1}^{N} \, I_{t}(k; \alpha^{2}) =
\prod\limits_{k=1}^{N} \, \left( \frac{\alpha^{2}(b-a)}{2N} \,
\int_{0}^{t} \, \sigma^{2} (u_{k};\tau) \, d \tau  + 1 \right)^{-
{1/2}} .
\end{equation}
In order to evaluate $ I_{1}(a,b;t) $, we take the logarithm of
\eqref{18}:
\begin{equation*}
\begin{array}{c}
\displaystyle\ln {\bf M} \, \left[ {I_{N}(k;\alpha^{2})}\right] =
\ln \prod\limits_{k=1}^{N} \, \left( \frac{\alpha^{2}(b-a)}{2N} \,
\int_{0}^{t} \, \sigma^{2} (u_{k};\tau) \, d \tau  + 1 \right)^{-
{1/2}} =
\\
= - \displaystyle\frac {1}{2} \, \sum \limits _{k=1}^{N} \, \ln
\left( \frac{\alpha^{2}(b-a)}{2N} \, \int_{0}^{t} \, \sigma^{2}
(u_{k};\tau) \, d \tau  + 1 \right) .
\end{array}
\end{equation*}
By using the series expansion of $ \ln{(x+1)} $ we obtain that:
\begin{equation*}
\begin{array}{c}
\displaystyle\ln{{\bf M} \, \left[ {\tilde I}_{1}(b;t) \right] } =
- \frac{1}{2} \, \sum \limits _{k=1}^{N} \, \left[
\frac{\alpha^{2}(b-a)}{2N} \, \int_{0}^{t} \, \sigma^{2}
(u_{k};\tau ) \, d \tau  - \right.
\\
\displaystyle\left.- \frac {1}{2}
\frac{\alpha^{4}(b-a)^{2}}{4N^{2}} \, \left( \int_{0}^{t} \,
\sigma^{2} (u_{k};\tau) \, d \tau \right)^{2} + O (N^{-3}) \right]
.
\end{array}
\end{equation*}
We calculate the limit for $N \to \infty$:
\begin{equation*}
\begin{array}{c}
\displaystyle\lim_{N \to \infty} \, \ln {\bf M} \, \left[
I_{N}(k,\alpha^{2}) \right] = - \frac {1}{2} \, \lim_{N
\rightarrow\infty} \, \sum \limits _{k=1}^{N} \, \left[
\frac{\alpha^{2}(b-a)}{2N} \, \int_{0}^{t} \, \sigma^{2}
(u_{k};\tau) \, d \tau  - \right.
\\
\displaystyle\left. - \frac {1}{2} \cdot
\frac{\alpha^{4}(b-a)^{2}}{4N^{2}} \, \left( \int_{0}^{t} \,
\sigma^{2} (u_{k};\tau) \, d \tau \right)^{2} + O (N^{-3}) \right]
=
\\
=  - \displaystyle\frac {1}{2} \, \lim_{N \to \infty} \, \sum
\limits _{k=1}^{N} \, \frac{\alpha^{2}(b-a)}{2N} \, \int_{0}^{t}
\, \sigma^{2} (u_{k};\tau) \, d \tau = - \frac{\alpha^{2}}{4} \,
\int_{a}^{b} \left( \int_{0}^{t} \, \sigma^{2} (u;\tau) \, d \tau
\right) \, du  .
\end{array}
\end{equation*}
By using the limit and passing to the anti-logarithm we prove that
Lemma \ref{lm2} holds
$$
{\bf M} \, \left[ \exp{ \left\{ - \frac{\alpha^{2}}{2}
\int_{a}^{b} \eta^{2} (u;t)\, du \right\} } \right] = \exp{
\left\{ - \frac{\alpha^{2}}{4} \int_{a}^{b} \left( \int_{0}^{t}
\sigma^{2}(u; \tau ) \, d \tau \right) \, du \right\} } . \ \ \
\square
$$

\subsection{Passing to the limiting process}

The field  $\{ x_{n}(l;t)$;  $ y_{n}(l;t) \} $ is defined by the
model assumption (\ref{3}) :
$$
 \eta (l_{s};t) = \int_{0}^{t} \sigma (l_{s}; \tau ) \, d
 w(\tau).
$$
We change the model assumption
$$
%\begin{equation}\label{111}
\displaystyle {\tilde \eta} (l_{s};t) = \left(\frac
{1}{2}\int_{0}^{t} \sigma^ {2} (l_{s}; \tau ) \, d \tau \right)
^{1/2}
%\end{equation}
$$
and consider the field
$ \{ {\hat x}_{n}(l;t) ; $ $ {\hat
y}_{n}(l;t) \} $
of the following form:
\begin{equation}\label{19}
\begin{array}{c}
\displaystyle {\hat x}_{n}(l;t) = \sum \limits_{s=1}^{n} \, a
(l_{s}) \cos \left[\sum \limits_{j=1}^{s} \, \triangle (w (l_{j}))
\, \left(\frac {1}{2}\int_{0}^{t} \sigma^{2} (l_{j}; \tau ) \, d
\tau \right) ^{1/2} \right] \ \Delta ,
\\
\displaystyle {\hat y}_{n}(l;t) = \sum \limits_{s=1}^{n} \, a
(l_{s}) \sin \left[\sum \limits_{j=1}^{s} \, \triangle (w (l_{j}))
\, \left(\frac {1}{2}\int_{0}^{t} \sigma^{2} (l_{j}; \tau ) \,
 d \tau
\right)^{1/2} \right] \ \Delta ,
\end{array}
\end{equation}
where $ \triangle (w (l_{j})) $ is the increment of the Wiener
process on the interval $ [l_{j};l_{j+1} ] $. This means that the
variable $t$ is not a random variable and from the analysis of the
process on the flow of the $\sigma$-algebras $ \displaystyle
\Im_{t}(n) \, \oplus \, \Im_{l} $ it is possible to pass to the
process defined on the flow of the $\sigma$-algebras $ \Im (l)$,
for all $t={\it const}$ . Averaging with respect to $t$ has
already been carried out. We observe that the fields $ \{
x_{n}(l;t); $ $ y_{n}(l;t) \} $ and $ \{ {\hat x}_{n}(l;t) ; $ $
{\hat y}_{n}(l;t) \} $ are defined on different spaces. We
consider the processes
$$
z_{2} (k;t) = \exp { \left \{ - {\it i} \, \sum \limits_{j=1}^{k}
\, \triangle w (l_{j}) \, \left(\frac {1}{2}\int_{0}^{t}
\sigma^{2} (l_{j}; \tau ) \, d \tau  \right) ^{1/2} \right \} } ,
$$
$$
z_{n,2} (l;t) = \sum \limits_{k=1}^{n} \, z_{2} (k;t) \, a(l_{k})
\cdot \Delta , \ \ \ \Delta = O (n^{-1}) .
$$
By considering the Euler representation, the components of the
fields (\ref{19})  will take the form:
$$
\begin{array}{c}
{\hat x}_{n}(l;t) = \displaystyle\frac {1}{2} \left ( z_{n,2}
(l;t) +
z_{n,2}^{*} (l;t) \right) , \\
{\hat y}_{n}(l;t) =  \displaystyle \frac {{\it i}}{2} \left (
z_{n,2} (l;t) - z_{n,2}^{*} (l;t) \right) .
\end{array}
$$
We construct now the characteristic functions: $ g_{n} ( \alpha ;
\beta ; t) $ for the field $ \left\{ x_{n}(l;t) ; y_{n}(l;t)
\right\} $ and $ {\hat g}_{n} ( \alpha ; \beta ; t) $ for the
field $ \left\{ {\hat x}_{n}(l;t) ; {\hat y}_{n}(l;t) \right\} $:
$$
g_{n} ( \alpha ; \beta; t) = {\bf M} \left[ \exp { \left \{ \frac
{i}{2} ( \alpha + i \, \beta) \,  z_{n,1} (l;t) + \frac {i}{2} (
\alpha - i \, \beta) \, z_{n,1}^{*} (l;t) \right \} } \right] ,
$$
$$
{\hat g}_{n} ( \alpha ; \beta; t) = {\bf M} \left[ \exp { \left \{
\frac {{\it i}}{2} ( \alpha + i \, \beta) \,  z_{n,2} (l;t) +
\frac {i}{2} ( \alpha - i \, \beta) \, z_{n,2}^{*} (l;t) \right\}
} \right] .
$$
For the continuation of the research the next lemma is necessary.

\begin{lem}\label{lm3}
{\it  Under the model assumptions for the random fields $ \left\{
x_{n}(l;t) ; y_{n}(l;t) \right\} $ and $ \left\{ {\hat x}_{n}(l;t)
; {\hat y}_{n}(l;t) \right\}  $ and for fixed integer
$ m$ there exists a number $%
n^{\prime}$ such that for all $n>n^{\prime}$ the following
relationships hold:
\begin{equation}\label{20}
{\bf M} [z_{,2}^{m}(l;t)] = {\bf M} [z_{,2}^{* \,m}(l;t)] =
 {\bf M} [z_{,1}^{m}(l;t)] =
  {\bf M} [z_{,1}^{* \, m}(l;t)] .
\end{equation}
}\end{lem}
 {\it Proof.} We calculate all increments in the series of the
equalities \eqref{20}  by having in mind the lemmas. We introduce
the following notation $\tilde{\eta}(\theta,t)=
\left(\displaystyle\frac{1}{2}\int_{0}^{t}\sigma^{2}(\theta,\tau)d\tau\right)^{1/2}$.
We have that
\begin{equation*}
\begin{array}{c}
{\bf M} [z_{,2}^{m}(l;t)]  = m !  \displaystyle\fint_{0}^{l} \, a
(u_{1}) \, d u_{1} \, {\bf M} \left[
 \exp {
\left\{ - i\, m \, \fint_{0}^{u_{1}} \tilde{\eta}(\theta_{1},t)\,
d w(\theta_{1}) \right\} } \right] \times
%\\
%\times \, \displaystyle\fint_{u_{1}}^{l} \, a (u_{2}) \, du_{2} \,
%{\bf M} \left[\exp { \left\{ - i\, (m-1) \, \fint_{u_{1}}^{u_{2}}
%\tilde{\eta}(\theta_{2},t)\, d w(\theta_{2}) \right\} }
%\right] \times \, \ldots \, \times
\\
\times\, \ldots \, \times \, \displaystyle\fint_{u_{m}}^{l} \, a
(u_{m}) \, du_{m}) \, {\bf M} \left[ \exp { \left \{ - i \,
\fint_{u_{m-1}}^{u_{m}} \tilde{\eta}(\theta_{m},t)\, d
w(\theta_{m}) \right\} } \right] =
\\
=m! \displaystyle\left(\fint_{0}^{l} \, a (u_{1}) \, d u_{1} \,
 \exp {
\left\{ - \frac{m^{2} }{2}\, \fint_{0}^{u_{1}}
\tilde{\eta}^{2}(\theta_{1},t)\, d\theta_{1} \right\} }
\right)\times
%\\
%\times \, \displaystyle\left(\fint_{u_{1}}^{l} \, a (u_{2}) \,
%du_{2} \, \exp { \left\{ - \frac{(m-1)^{2}}{2} \,
%\fint_{u_{1}}^{u_{2}} \tilde{\eta}^{2}(\theta_{2},t)\, d
%\theta_{2} \right\} }
% \right)\times \, \ldots \, \times
\\
\times \,  \ldots \, \times\,\displaystyle\left(\fint_{u_{m}}^{l}
\, a (u_{m}) \, du_{m}) \,  \exp { \left \{ - \frac{1}{2} \,
\fint_{u_{m-1}}^{u_{m}} \tilde{\eta}^{2}(\theta_{m},t)\, d
\theta_{m} \right\} } \right)=
\\
=m! \displaystyle\left(\fint_{0}^{l} \, a (u_{1}) \, d u_{1} \,
 \exp {
\left\{ - \frac{m^{2} }{2}\, \fint_{0}^{u_{1}}
\left(\displaystyle\frac{1}{2}\int_{0}^{t}\sigma^{2}(\theta_{1},\tau)d\tau\right)\,
d\theta_{1} \right\} }  \right)\times \,\ldots \, \times
%\\
%\times \, \displaystyle\left(\fint_{u_{1}}^{l} \, a (u_{2}) \,
%du_{2} \, \exp { \left\{ - \frac{(m-1)^{2}}{2} \,
%\fint_{u_{1}}^{u_{2}}
%\left(\displaystyle\frac{1}{2}\int_{0}^{t}\sigma^{2}(\theta_{2},\tau)d\tau\right)\,
%d \theta_{2} \right\} }
% \right)\times \, \ldots \, \times
\\
 \times\, \displaystyle\left(\fint_{u_{m}}^{l}
\, a (u_{m}) \, du_{m}) \,  \exp { \left \{ - \frac{1}{2} \,
\fint_{u_{m-1}}^{u_{m}}
\left(\displaystyle\frac{1}{2}\int_{0}^{t}\sigma^{2}(\theta_{m},\tau)d\tau\right)\,
d \theta_{m} \right\} } \right)={\bf M} [z_{,1}^{m}(l;t)]
\end{array}
\end{equation*}
Then, in force of \eqref{7}  and \eqref{9}  we have that:
$$
{\bf M} [z_{,1}^{m}(l;t)] = {\bf M} [z_{,1}^{* \,m}(l;t)]. \ \
$$
In this way we obtain the confirmation of Lemma.\ \ \ $\square$

\begin{lem}\label{lm4}
{\it The characteristic functions of fields $ \{ x_{n}(l;t)$; $
y_{n}(l;t) \} $ and $ \{ {\hat x}_{n}(l;t)$; $ {\hat y}_{n}(l;t)
\} $ for $ n \to \infty $ coincide for all $l$ and $t$.}
\end{lem}

{\it Proof.} The proof is based on the coincidence of the
representations for the characteristic functions $ g_{n}(\alpha ,
\beta , t) $ and  $ {\hat g}_{n}(\alpha , \beta , t) $ by means of
the Maclaurin expansion (inside the mean) with respect to $
z_{,1}(l;t) $ and $ z_{,1}^{*}(l;t) $, and $ z_{,2}(l;t) $, $
z^{*}_{,2}(l;t)$ respectively and also on the conclusions of Lemma
\ref{lm3}. \ \ \ $\square$

 Lemma \ref{lm4} permits us to pass to
the study of the limit behavior of the field $ \{ {\hat
x}_{n}(l;t)$; $ {\hat y}_{n}(l;t) \}$ for $n \to \infty$
exclusively.

\begin{thm}\label{t1} Let us assume that for the field $ \{
x_{n}(l;t); y_{n}(l;t) \} $ the model assumptions {\rm(\ref{19})}
are satisfied:
\begin{equation*}
\begin{array}{c}
\displaystyle {\hat x}_{n}(l;t) = \sum \limits_{s=1}^{n} \, a
(l_{s}) \cos \left[\sum \limits_{j=1}^{s} \, \triangle(w (l_{j}))
\, \left(\frac {1}{2}\int_{0}^{t} \sigma^{2} (l_{j}; \tau ) \, d
\tau \right) ^{1/2} \right] \ \Delta ,
\\
\displaystyle {\hat y}_{n}(l;t) = \sum \limits_{s=1}^{n} \, a
(l_{s}) \sin \left[\sum \limits_{j=1}^{s} \,\triangle (w (l_{j}))
\, \left(\frac {1}{2}\int_{0}^{t} \sigma^{2} (l_{j}; \tau ) \,
 d \tau
\right)^{1/2} \right] \ \Delta ,
\end{array}
\end{equation*}
and assume that the field $ \{ x(l;t); y(l;t) \} $ is
defined in the following way
\begin{equation}\label{21}
\begin{array}{c}
x(l;t)=\displaystyle \int_{0}^{l} a(u) \cos \left[\int_{0}^{u}
\frac{1}{2}\left(\int_{0}^{t}\sigma^{2}(\theta,\tau)d\tau\right)dw(\theta)
\right]du,\\
y(l;t)=\displaystyle \int_{0}^{l} a(u) \sin \left[\int_{0}^{u}
\frac{1}{2}\left(\int_{0}^{t}\sigma^{2}(\theta,\tau)d\tau\right)dw(\theta)
\right]du.
\end{array}
\end{equation}
Under these conditions the characteristic functions of the
processes $ \{ x(l;t); y(l;t) \} $ and $ \{ x_{n}(l;t); y_{n}(l;t)
\} $ coincide.
\end{thm}

{\it Proof.} The comparison of the characteristic functions for $
\{ {\hat x}_{n}(l;t), {\hat y}_{n}(l;t) \} $ and $ \{ x(l;t), $ $
y(l;t) \} $ for all values of $t \in [0,T],$ for $n \to \infty$
leads to the proof of the theorem.\ \ \ $\square$

\begin{thm}\label{t2} The stochastic process $ \{ x(l;t); y(l;t) \} $
is the solution to the Cauchy problem for the Ito stochastic
differential equations:
\begin{equation}\label{22}
\begin{array}{c}
\displaystyle d_{l} p(l;t) = \left[ p(l;t) \frac{\partial
}{\partial l} \ln  a(l) - \frac{p(l;t)}{4}  \int \limits _{0}^{t}
\sigma ^{2}(l; \tau )d \tau \right]  dl %- \\
 - \displaystyle\left( \frac{1}{2}
\int\limits _{0}^{t} \sigma ^{2}(l; \tau ) d \tau \right)^{0,5}
q(l;t)  dw(l),
\\
\displaystyle d_{l} q(l;t) = \left[ q(l;t)  \frac{\partial
}{\partial l} \ln  a(l)  - \frac{q(l;t)}{4}  \int \limits
_{0}^{t} \sigma ^{2}(l; \tau ) d \tau \right] dl %+ \\
+ \displaystyle\left( \frac{1}{2} \int \limits _{0}^{t} \sigma
^{2}(l; \tau ) d \tau \right)^{0,5} p(l;t)  dw(l),
\end{array}
\end{equation}
$$
d_{l}x(l;t) = q(l;t)\, dl ,  \ \ \ \ d_{l}y(l;t) = p(l;t)\, dl ,
$$
satisfying the boundary conditions
$$
x(0;t)=0, \ \ y(0;t)=0, \ \ p(0;t)=a(0), \ \ q(0;t)=0.
$$
\end{thm}
{\it Proof.} We differentiate $ x(l;t) $ and $ y(l;t) $ in
\eqref{21} with respect to $l$:
\begin{equation}\label{23}
\displaystyle \frac{\partial x(l;t)}{\partial l}=a(l) \sin
\left[\int_{0}^{l}
\frac{1}{2}\left(\int_{0}^{t}\sigma^{2}(\theta,\tau)d\tau\right)dw(\theta)
\right]=q(l;t),
\end{equation}
\begin{equation}\label{24}
\displaystyle \frac{\partial y(l;t)}{\partial l}=-a(l) \cos
\left[\int_{0}^{l}
\frac{1}{2}\left(\int_{0}^{t}\sigma^{2}(\theta,\tau)d\tau\right)dw(\theta)
\right]=p(l;t).
\end{equation}
The obtained expressions are now differentiated by Ito formula
with respect to the variable $l$:
$$
\begin{array}{c}
 d_{l}\left(\displaystyle \frac{\partial x(l;t)}{\partial
l}\right)=\displaystyle\frac{1}{a(l)}\frac{\partial a(l)}{\partial
l}\,a(l)\cos \left[\int_{0}^{l}
\frac{1}{2}\left(\int_{0}^{t}\sigma^{2}(\theta,\tau)d\tau\right)dw(\theta)
\right]dl-
\\
-\,a(l) \sin \left[\displaystyle\int_{0}^{l}\displaystyle
\frac{1}{2}\left(\int_{0}^{t}\sigma^{2}(\theta,\tau)d\tau\right)dw(\theta)
\right]\cdot \displaystyle\frac{1}{2}
\left(\int_{0}^{t}\sigma^{2}(\theta,\tau)d\tau\right)dw(l)-
\\
-\,\displaystyle\frac{1}{2}\cos \left[\int_{0}^{l}
\frac{1}{2}\left(\int_{0}^{t}\sigma^{2}(\theta,\tau)d\tau\right)dw(\theta)
\right]\cdot\left(\displaystyle\frac{1}{2}
\int_{0}^{t}\sigma^{2}(\theta,\tau)d\tau\right)^{2}dl.
\end{array}
$$
Taking into account \eqref{23}  and \eqref{24}  we obtain the last
equation of system \eqref{22}. In a similar way we get the second
expression of the system. The functions $x(l;t)$, \ $y(l;t)$, \
$p(l;t)$, \ $q(l;t)$ defined by \eqref{21}, \eqref{23}  and
\eqref{24} satisfy the given initial conditions.\ \ \ $\square$

Within the framework of the given formulation ($L$ constant) we
have found $ F_{t}(x;y;L) $ for different values of $t$.

\begin{thm}\label{t3}  The distribution function of the process $ \{
x(l;t); y(l;t) \} $ can be obtained by integrating with respect to
the variables $p$ and $q$ the Kolmogorov equation of the system
 {\rm(\ref{22})}.
\end{thm}

{\it Proof.} After the enlargement of the space obtained by
introducing the new variables $p$ and $q$, the compound process $
\{ x(l;t); y(l;t); p(l;t);q(l;t) \} $ becomes a Markov process.
This means that it is possible to obtain a Kolmogorov equation for
the density function $\rho (x,y,p,q,l,t)$ and then by integrating
with respect to $p$ and $q$ infer the density function of the
distribution $\rho (x,y,l,t)$ for all $l$ and $t$.\ \ \ $\square$

\begin{thm}\label{t4} The distribution function of the original
process $ \{ x_{n}(l;t); y_{n}(l;t) \} $ under the model
conditions {\rm(\ref{4})}  coincide with the distribution function
of the Markov process $ \{ x(l;t); y(l;t) \} $ {\rm(\ref{21})}.
\end{thm}

{\it Proof.} The proof is based on the conclusions of Theorem
\ref{t1} and Theorem \ref{t2}.\ \ \ $\square$
\newline

\begin{rem} The character of the analysis doesn't
substantially changes when, for example $a=a(l,t)$ (vibrating
chain), $ \sigma (l;t) $ is a non anticipating measurable random
function with respect to independent flows of $\sigma$-algebras
governed by independent Wiener processes $ w(l) $ and $ w(t)$.
\end{rem}

In this way we arrive at a coherent representation of
distribution: the parameter $t$ defines also the structure of the
chain.

\end{document}